\documentclass[12pt]{article}

\usepackage{amsmath}
\usepackage{amssymb}
\usepackage{amsthm}
\usepackage{xypic}
\usepackage{url}
\usepackage{dsfont}

\begin{document}
\baselineskip=16pt
\textheight=9.3in
\parindent=0pt 
\def\sk {\hskip .5cm}
\def\skv {\vskip .12cm}

\large\bf\centerline{A Commutative Alternative to Fractional Calculus on $k$-Differentiable Functions}\rm\normalsize
\vskip .75cm
\centerline{Matthew Parker}
\vskip 1.5cm

\skv
\large{\bf Abstract}\normalsize
\vskip .1cm
\qquad In [1], an operator was introduced which acts parallel to the Riemann-Liouville differintegral $_{a}D_{x}^{k}f(x) = \frac{1}{\Gamma(-k)}\int_{a}^{x}f(t)(x-t)^{-\alpha-1}dt$ [2] on a transformation of the space of real analytic functions $C^{\omega}$, denoted by $\mathbb{R}_{\omega}$, and commutes with itself. This paper aims to extend the technique - and its defining characteristic, commutativity - to all real continuous functions, up to the degree to which they are differentiable.
\vskip 1cm

\skv
\large{\bf Convention}\normalsize
\vskip .1cm
\qquad We will let $C(\mathbb{R})$ denote the set of continuous functions $\mathbb{R} \to \mathbb{R}$, and $C^{k}$ denote real-valued functions $f$ for which $k$ is the largest integer such $\frac{d^{k}}{dx^{k}}f$ exists and is continuous everywhere. We write $_{a}I_{x}^{k}f(x)$ for the Riemann-Liouville integral of degree $k$, and the term "derivative" shall refer to the unrestricted R-L derivative; that is, for all degrees $p \in \mathbb{R}\cup\{\infty, \omega\}$. We define, for $f \in C^{k}$ and $a \in \mathbb{R}, \,_{a}T(f) = \sum_{i=0}^{k}\frac{[\frac{d^{i}}{dx^{i}}f](a)}{i!}(x-a)^{i}$. We adopt the definitions of $\mathbb{R}_{\omega}, D^{k}, R,$ and $R^{-1}$ from [1].
\vskip 1cm

\skv
\large{\bf Introductory Discussion}\normalsize
\vskip .1cm
\qquad From [1], we have $\mathbb{R}_{\omega} = \{\sigma : \mathbb{R} \to \mathbb{R}$ such that $\sum_{i=0}^{\infty}\frac{\sigma(i)}{\Gamma(i+1)}(x-a)^{i}$ converges on some subset of $\mathbb{R}$\}, and defined the operator $D^{k}$ to act on $\mathbb{R}_{\omega}$ as a shift operator; that is, $[D^{k}\sigma](i) = \sigma(i-k)$. This had the convenient property that $\sum_{i=0}^{\infty}\frac{D^{k}\sigma(i)}{\Gamma(i+1)}(x-a)^{i} = \frac{d^{k}}{dx^{k}}\sum_{i=0}^{\infty}\frac{\sigma(i)}{\Gamma(i+1)}(x-a)^{i}$. Maps between the space of analytic functions and $\mathbb{R}_{\omega}$ allowed this shift operator to work parallel to the Riemann-Liouville derivative. If we consider continuous functions as our originating space, however, the ability to condense all information about a function into a countable collection (e.g, $f \leftrightarrow \{\frac{d^{i}}{dx^{i}}f\vert_{x=a} : i \in \mathbb{N}$\}) disappears, and more information must be retained.
\vskip 1cm

\skv
\large{\bf Creating a Commutative Alternative}\normalsize
\vskip .1cm
\qquad For each $f \in C^{k}$ and nonnegative $p, q$ with $p \geq q$, we know from [2] and [4] that $_{a}D_{x}^{-q}[_{a}D_{x}^{p}f] = \,_{a}D_{x}^{p-q}f$ modulo a sum of the form $\sum_{i=0}^{\lfloor\,\vert p\vert\,\rfloor}c_{i}x^{i-q}$. As differentiation forms a semigroup [3], this is exactly the failure of the commutativity of the derivative. To obtain a commutative operator working parallel to the derivative, we must specify exactly what this polynomial is. \\*

\qquad If $a$ is a point where $f$ is not $k+1$-differentiable, and we let $p, q \to k, k+q-p$ keeping $p-q$ constant, using integration by parts $k$ times we see in the limit, $_{a}D_{x}^{-q}[_{a}D_{x}^{p}f] = \,_{a}D_{x}^{p-q}f - \,_{a}D_{x}^{p-q}[_{a}T(f)]$ [4]. Observing $_{a}D_{x}^{p-q}[_{a}T(f)] = \,_{a}D_{x}^{p}RD^{-q}\sigma$ for some $\sigma \in \mathbb{R}_{\omega}$, and noting from [1] we have a method for differentiating commutatively with all such $\sigma \in \mathbb{R}_{\omega}$, this suggests we may begin to form a commutative view of differentiation of $f$ if we consider it to be comprised of two parts; a function $f_{a} \in C^{k}$ with $_{a}T(f_{a}) = 0$ whenever $f$ is not differentiable at $a$, and some $\sigma \in \mathbb{R}_{\omega}$. \\*

\qquad We required $f$ not be $k+1$-differentiable at the point $a$ because it allows us to construct the equation $_{a}D_{x}^{-q}[_{a}D_{x}^{p}f] = _{a}D_{x}^{p-q}f - _{a}D_{x}^{p-q}[_{a}T(f)]$, which in turn will allow us to separate out the polynomial in question. Clearly, if $f$ were $n$-differentiable at the point $a$ for some $n \geq k+1$, then $_{a}T(f)$ would be a polynomial of degree $n$, and $_{a}D_{x}^{p-q}f - _{a}D_{x}^{p-q}[_{a}T(f)]$ = $_{a}D_{x}^{-q}[_{a}D_{x}^{p}f] - P_{n-k}(x)$, where $P_{n-k}(x)$ is some polynomial of degree $n-k$. \\*

\qquad The equation $_{a}D_{x}^{-q}[_{a}D_{x}^{p}f] = _{a}D_{x}^{p-q}f - _{a}D_{x}^{p-q}T(f)$ shows us exactly how to define the corresponding pair ($f_{a}, \sigma$), and the proper actions on them:
\vskip .8cm

{\bf Definition 1.a:} For $f \in C^{k}$, if $a$ is a point where $f$ is not $k+1$-differentiable, define $r^{-1}f$ be the pair ($f_{a}, \sigma$) such that $f_{a} = \,_{a}D_{x}^{-k}[_{a}D_{x}^{k}f] \in C^{k}$, and $\sigma = R^{-1}[_{a}T(f)] \in \mathbb{R}_{\omega}$. 
\vskip .8cm

{\bf Definition 1.b:}For a pair ($f_{a}, \sigma)$ such that $_{a}D_{x}^{-k}[_{a}D_{x}^{k}f_{a}] = f_{a}$ where $f_{a} \in C^{k}$ and $\sigma \in \mathbb{R}_{\omega}$, define $r(f_{a}, \sigma) = f_{a} + R\sigma$. 
\vskip .8cm

{\bf Definition 2:} Denote the space $rC^{k} = \{(f_{a}, \sigma) : f_{a} = \,_{a}D_{x}^{-k}[_{a}D_{x}^{k}f]$ for some $f \in C^{k}, \sigma \in \mathbb{R}_{\omega}, \sigma(i) = 0$ for $i \textgreater\, k\}$ by $C_{k}$.
\vskip .8cm

{\bf Definition 3:} For $f \in C^{k}$ and $p \,\,\textgreater\, k-1$, if ($f_{a}, \sigma$) = $r^{-1}f$, define $D^{p}(f_{a}, \sigma) = (_{a}D_{x}^{p}f_{a}, D^{p}\sigma)$.
\vskip 3cm

\skv
\large{\bf Properties of $D^{k}, C_{k}, r$, and $r^{-1}$}\normalsize
\vskip .1cm

\qquad We claim the operator $D^{p}$ acts parallel to the Riemann-Liouville derivative in the sense that the following diagram commutes, for $f \in C^{k}$:
\begin{displaymath}
    \xymatrix{f \ar@/_3pc/[dd]_{_{a}D_{x}^{p+q}} \ar[d]^{_{a}D_{x}^{p}} \ar[rr]^{r^{-1}} & & (f_{a}, \sigma) \ar[d]^{D^{p}} \\
               _{a}D_{x}^{p}f & & \ar[ll]^{r} D^{p}(f_{a}, \sigma)\ar[d]^{D^{q}} \\
               _{a}D_{x}^{p+q}f & &\ar[ll]^{r} D^{q}D^{p}(f_{a}, \sigma)  }
\end{displaymath}\\*

\qquad Let $f \in C^{k}$, and $(f_{a}, \sigma) = r^{-1}f$. Let $p, q$ be such that $p, q, p+q\,\, \textless\, k+1$. Then
\begin{align*}
rD^{q}D^{p}r^{-1}f \,\, & =  \,rD^{q}D^{p}(f_{a}, \sigma) \\
\\
 & = \,rD^{q}(_{a}D_{x}^{p}f_{a}, D^{p}\sigma) \\
 \\
 & = \,r(_{a}D_{x}^{q}[_{a}D_{x}^{p}f_{a}], D^{q}D^{p}\sigma) \\
 \\
 & =\, r(_{a}D_{x}^{q+p}f_{a}, D^{q+p}\sigma) \\
 \\
 & = \, _{a}D_{x}^{q+p}[_{a}D_{x}^{-k}[_{a}D_{x}^{k}f_{a}]] + \sum_{i=0}^{k-q-p}\frac{\sigma(i)}{i!}(x-a)^{i} \\
 \\
 & = \,_{a}D_{x}^{p+q}f +\, _{a}D_{x}^{p+q}T(f) \\
 \\
 & = \, _{a}D_{x}^{p+q}f 
\end{align*}

and we see the diagram commutes. \\*

\qquad We next show the operator $D^{k}$ commutes with itself. Let $(f_{a}, \sigma) \in C_{k}$ and let $p, q$ be such that $p, q, p+q\,\, \textless\, k+1$. Then \\*

\centerline{$D^{q}D^{p}(f_{a}, \sigma) \,\, = \,\, (D^{q}_{a}D_{x}^{p}f_{a}, D^{p+q}\sigma)$}
\vskip .5cm
\qquad Noting $f_{a} = \,_{a}D_{x}^{-k}[_{a}D_{x}^{k}f]$ for some $f \in C^{k}$, we see
\vskip 2cm
\begin{align*}
D^{q}D^{p}(f_{a}, \sigma) \, & =\, (_{a}D_{x}^{q}[_{a}D_{x}^{p}[_{a}D_{x}^{-k}[_{a}D_{x}^{k}f_{a}]]], \, D^{q}D^{p}\sigma) \\
\\
& =\,(_{a}D_{x}^{p+q}[_{a}D_{x}^{-k}[_{a}D_{x}^{k}f_{a}], \, D^{p+q}\sigma)\\
\\
& =\,(_{a}D_{x}^{q+p}[_{a}D_{x}^{-k}[_{a}D_{x}^{k}f_{a}]], \, D^{q+p}\sigma) \qquad\qquad \textnormal{by definition of}\, f_{a}\\
\\
& =\,(_{a}D_{x}^{p}[_{a}D_{x}^{q}[_{a}D_{x}^{-k}[_{a}D_{x}^{k}f_{a}]]], \, D^{p}D^{q}\sigma)  \,\,\qquad \textnormal{since}\,\, p, q, p+q \textless\, k+1 \qquad [2]  \\
\\
& =\,(_{a}D_{x}^{p}[_{a}D_{x}^{q}f_{a}], \, D^{q}D^{p}\sigma) \\
\\
& =\,D^{p}D^{q}(f_{a}, \sigma)
\end{align*}
\vskip 1cm

\skv
\large{\bf Unifying the Spaces $C_{k}$}\normalsize
\vskip .1cm
\qquad We have thus far developed a collection of spaces $C_{k}$ and operators $D^{p}$ on those spaces. We now attempt to unify these spaces in an intuitive way, put a vector space structure on the result, and define an operator $D^{p}$ which acts parallel to the Riemann-Liouville derivative and commutes with itself. \\*

{\bf Definition 4:} Define $C_{diff} = \{(f_{a}, \sigma) : f_{a}$ continuous, $f_{a}$ not $k+1$-differentiable at $a$ if $f \in C^{k}, \,\,_{a}D_{x}^{-k}[_{a}D_{x}^{k}f_{a}] = f_{a}$ if $f_{a} \in C^{k}$, and $\sigma \in \mathbb{R}_{\omega}$\}. \\*

{\bf Definition 5.a:} For $(f_{a}, \sigma), (g_{a}, \rho) \in C_{diff}$, define $(f_{a}, \sigma) + (g_{a}, \rho) = $

($f_{a} + g_{a}, \sigma + \rho$). \\*

{\bf Definition 5.b:} For $r \in \mathbb{R}$ and $(f_{a}, \sigma) \in C_{diff}$, define $r(f_{a}, \sigma) = (rf_{a}, r\sigma)$.\\*

\qquad Definitions 5.a, 5.b turn $C_{diff}$ into a vector space with identity (0, 0) since the sum and scalar multiple of continuous functions is continuous, and $\mathbb{R}_{\omega}$ is closed under addition and scalar multiplication. \\*

{\bf Definition 6:} For ($f_{a}, \sigma) \in C_{diff}$, if $f_{a} \in C_{k}$, define $D^{p}(f_{a}, \sigma) = (_{a}D_{x}^{p}f_{a}, D^{p}\sigma)$ for all $p \,\textless\, k-1$. \\*

\qquad From the discussions in the previous section, we see immediately $D^{p}$ commutes with itself.

\vskip 3cm

\skv
\large{\bf Final Constructions and Properties}\normalsize
\vskip .1cm
 \qquad We now attempt to relate $C_{diff}$ to $C(\mathbb{R})$. We will do this in much the same way we related $C^{k}$ to $C_{k}$. In fact, the method carries through almost exactly. \\*
 
 {\bf Definition 7.a:} For $f_{a} \in C(\mathbb{R}),$ if $f \in C^{k}$ for some $k \in \mathbb{R}\cup\{\infty, \omega\}$, and if $f$ is not $k+1$ - differentiable at $a$ if $k \in \mathbb{R}$, then define
 \vskip .1cm
\centerline{ $R^{-1}f = (_{a}D_{x}^{-k}[_{a}D_{x}^{k}f],\,\, R^{-1}[_{a}T(f)]) \in C_{diff}$}
\vskip .5cm
{\bf Definition 7.b:} For ($f_{a}, \sigma) \in C_{diff}$, define $R(f_{a}, \sigma) = f_{a} + R\sigma$. \\*

\qquad It is immediate the maps $R, R^{-1}$ are homomorphisms, and $R$ is injective. \\*

{\bf Definition 8:} For ($f_{a}, \sigma) \in C_{diff}$ where $f_{a} \in C^{k}$, define $D^{p}(f_{a}, \sigma) = (_{a}D_{x}^{p}f_{a}, D^{p}\sigma)$ whenever $p \,\,\textless\, k+1$. \\*

\qquad From the discussions of the earlier sections, clearly $D^{p}$ and $D^{q}$ commute on ($f_{a}, \sigma) \in C_{diff}$ whenever $p, q, p+q \,\,\textless\, k+1$ if $f_{a} \in C^{k}$. We now need only establish $_{a}D_{x}^{k}f = RD^{k}R^{-1}f$ to show equivalence of $D^{p}$ on $C_{diff}$ to $_{a}D_{x}^{p}$ on $C(\mathbb{R}$). \\*

\qquad Let $f \in C^{k}$ for some $k \in \mathbb{R}\cup\{\infty, \omega\}$, let $a$ be a point at which $f$ is not $k+1$-differentiable if $k \in \mathbb{R}$, and $p \,\, \textless\, k+1$. Then

\vskip .1cm
\begin{align*}
RD^{k}R^{-1}f &= \, RD^{p}(_{a}D_{x}^{-k}[_{a}D_{x}^{k}f], R^{-1}[_{a}T(f)]) \\
\\ 
& = \, R(_{a}D_{x}^{p}[_{a}D_{x}^{-k}[_{a}D_{x}^{k}f], D^{p}R^{-1}[_{a}T(f)]  \\
\\ 
& = \, _{a}D_{x}^{p}[_{a}D_{x}^{-k}[_{a}D_{x}^{k}f]] + RD^{p}R^{-1}[_{a}T(f)]  \\
\\ 
& = \, _{a}D_{x}^{p}[_{a}D_{x}^{-k}[_{a}D_{x}^{k}f] + _{a}T(f)]  \\
\\ 
& = \, _{a}D_{x}^{p}f
\end{align*}

\vskip 5cm

\skv
\large{\bf Remark}\normalsize
\vskip .1cm
\qquad It is worth noting that the entirety of the theory presented above carries through unchanged for $f \in C^{k}[a - \epsilon, b]$ for some $\epsilon \,\,\textgreater\, 0$ on the interval ($a - \epsilon, b)$, as the functions $f, _{a}D_{x}^{-k}[_{a}D_{x}^{k}f]$, and $_{a}T(f)$ will exist and satisfy the requisite equations on ($a-\epsilon, b)$.
\vskip 1cm

\skv
\large{\bf Conclusion}\normalsize
\vskip .1cm
\qquad In conclusion, we have created spaces \{$C_{k} : k \in \mathbb{R}\cup\{\infty, \omega$\} which contain images of the \{$C^{k}$\}, collected these together into a space $C_{diff}$ which contains a subset isomorphic to $C(\mathbb{R})$, and created an operator $D^{k}$ on $C_{diff}$ which commutes with itself, and acts on $C_{diff}$ such that the following diagram commutes;

\begin{displaymath}
    \xymatrix{ C(\mathbb{R}) \ar[d]_{_{a}D_{x}^{p}} \ar[r]^{R^{-1}} & C_{diff} \ar[d]^{D^{p}} \\
               C(\mathbb{R}) & C_{diff} \ar[l]^{R} }
\end{displaymath} \\*
\qquad Moreover, we showed the operator $D^{p}$ operates on $C_{diff}$ equivalently to $_{a}D_{x}^{p}$ on $C(\mathbb{R})$ up to the degree to which the function under consideration can be differentiated. This completes the paper.

\end{document}